\documentclass[12pt]{article}
\usepackage{amsmath,amssymb,amsthm}
\usepackage[mathscr]{eucal}
\usepackage[hypertex]{hyperref}

\newcommand{\arr}{\longrightarrow}
\newcommand{\B}{\mathcal{B}}
\newcommand{\pcf}{p.c.f.}
\newcommand{\xs}{X^*}
\newcommand{\xmo}{X^{-\omega}}

\newcommand{\lims}[1][G]{\mathscr{J}_{#1}}
\newcommand{\nuke}{\mathcal{N}}
\newcommand{\si}{\mathsf{s}}
\newcommand{\til}{\mathscr{T}}
\newcommand{\emp}{\emptyset}

\title{Post-critically finite self-similar groups}
\author{Ievgen Bondarenko, Volodymyr Nekrashevych\thanks{The second author acknowledges the support of Swiss National
Science Foundation and Alexander von Humboldt Foundation}}
\date{October 20, 2003}

\newtheorem{theorem}{Theorem}[section]
\newtheorem{proposition}[theorem]{Proposition}
\newtheorem{corollary}[theorem]{Corollary}

\theoremstyle{definition}
\newtheorem{defi}{Definition}[section]

\begin{document}
\maketitle

\textit{Dedicated to R.~I.~Grigorchuk on the occasion of his 50th
birthday}

\begin{abstract}
We describe in terms of automata theory the automatic actions with
post-critically finite limit space. We prove that these actions
are precisely the actions by bounded automata and that any
self-similar action by bounded automata is contracting.
\end{abstract}

\section{Introduction}

The aim of this paper is to show a connection between two notions,
which have appeared in rather different fields of mathematics. One
is the notion of a post-critically finite self-similar set (other
related terms are: ``nested fractal'' or ``finitely ramified
fractal''). It appeared during the study of harmonic functions and
Brownian motion on fractals. The class of post-critically finite
fractals is a convenient setup for such studies. See the
papers~\cite{kigami:anal_fract,kigami:lapl94,lind:Brown_motion,sabot:diffusion}
for the definition of a post-critically finite self-similar sets
and for applications of this notion to harmonic analysis on
fractals.

The second notion appeared during the study of groups generated by
finite automata (or, equivalently, groups acting on rooted trees).
Many interesting examples of such groups where found (like the
Grigorchuk group~\cite{gri:BurnsProb}, groups defined by
Aleshin~\cite{al:burn_en}, Sushchansky~\cite{susz:burn_en},
Gupta-Sidki groups~\cite{gupta-sidki_group} and many others), and
these particular examples where generalized to different classes
of groups acting on rooted trees: branch
groups~\cite{gri:justinf}, self-similar (state closed)
groups~\cite{fractal_gr_sets,sidki_monogr},
GGS-groups~\cite{baumslag:combgrth}, AT
groups~\cite{Merzl:periodgr,rozhkov:PhD}, spinal
groups~\cite{branch}.

S.~Sidki has defined in his work~\cite{sid:cycl} a series of
subgroups of the group of finite automata, using the notions of
activity growth and circuit structure. In particular, he has
defined the notion of a bounded automaton. The set of all
automorphisms of the regular rooted tree, which are defined by
bounded automata is a group. It is interesting that most of the
known interesting examples of groups acting on rooted trees (in
particular, all the above mentioned examples) are subgroups of the
group of bounded automata. Also every finitely automatic
GGS-group, AT-group or spinal group is a subgroup of the group of
bounded automata.

We prove in our paper that a self-similar (state closed) group is
a subgroup of the group of all bounded automata if and only if its
\emph{limit space} is a post-critically finite self-similar space.
The limit space of a self-similar group was defined
in~\cite{nekr:limsp} (see also~\cite{fractal_gr_sets}). This
establishes the mentioned above connection between the harmonic
analysis on fractals and group theory.

The structure of the paper is the following. Section
``Self-similar groups'' is a review of the basic definitions of
the theory of self-similar groups of automata. We define the
notions of self-similar groups, automata, Moore diagrams,
contracting groups, nucleus of a contracting group and establish
notations.

Third section ``Limit spaces'' gives the definition and the basic
properties of the limit space of a contracting self-similar group
as a quotient of the space of infinite sequences. We also discuss
the notion of \emph{tiles} of a limit space (the images of the
cylindric sets of the space of sequences).

The main results of the section ``Post-critically finite limit
space'' are Corollary~\ref{cor:pcflim}, giving a criterion when
the limit space of a self-similar group action is post-critically
finite and Proposition~\ref{pr:1dim}, stating that a
post-critically finite limit space is 1-dimensional.

The last section ``Automata with bounded cyclic structure'' is the
main part of the article. We prove Theorem~\ref{th:main}, which
says that every self-similar subgroup of the group $\B$ of bounded
automata is contracting and that a contracting group has a
post-critically finite limit space if and only if it is a subgroup
of $\B$.

\section{Self-similar groups}

We review in this section the basic definitions and theorems
concerning self-similar groups. For a more detailed account,
see~\cite{fractal_gr_sets}.

Let $X$ be a finite set, which will be called \emph{alphabet}. By
$\xs$ we denote the set of all finite words $x_1x_2\ldots x_n$
over the alphabet $X$, including the empty word $\emp$.

\begin{defi}
A faithful action of a group $G$ on the set $\xs$ is
\emph{self-similar} (or \emph{state closed}) if for every $g\in G$
and for every $x\in X$ there exist $h\in G$ and $y\in X$ such that
\[
g(xw)=yh(w)
\]
for all $w\in\xs$.
\end{defi}

We will write formally
\begin{equation}
\label{eq:ssdot} g\cdot x=y\cdot h,
\end{equation}
if for every $w\in\xs$ we have $g(xw)=yh(w)$. If one identifies
every letter $x\in X$ with the map $w\mapsto xw:\xs\arr\xs$, then
equation~\eqref{eq:ssdot} will become a correct equality of two
transformations.

The notion of a self-similar action is closely related with the
notion of an \emph{automaton}.

\begin{defi}
An \emph{automaton} $\mathcal{A}$ over the alphabet $X$ is a tuple
$\langle Q, \pi, \lambda\rangle$, where $Q$ is a set (the set of
\emph{internal states} of the automaton), and $\pi:Q\times X\arr
Q$ and $\lambda:Q\times X\arr Q$ are maps (the \emph{transition}
and the \emph{output functions}, respectively).

An automaton is \emph{finite} if its set of states $Q$ is finite.
A subset $Q'\subset Q$ is called \emph{sub-automaton} if for all
$q\in Q'$ and $x\in X$ we have $\pi(q, x)\in Q'$. If $Q'$ is a
sub-automaton, then we identify it with the automaton $\langle Q',
\pi|_{Q'\times X}, \lambda|_{Q'\times X}\rangle$.
\end{defi}

For every state $q\in Q$ and $x\in X$ we also write formally
\begin{equation}
\label{eq:automatdot} q\cdot x=y\cdot p,
\end{equation}
where $y=\lambda(q, x)$ and $p=\pi(q, x)$.

We will also often use in our paper another notation for the
functions $\pi$ and $\lambda$:
\[
\pi(q, x)=q|_x, \qquad \lambda(q, x)=q(x).
\]

The transition and output functions are naturally extended to
functions $\pi:Q\times\xs\arr Q$ and $\lambda:Q\times\xs\arr\xs$
by the formulae:
\[
\pi(q, xv)=\pi\left(\pi(q, x), v\right),\qquad\lambda(q,
xv)=\lambda(q, x)\lambda\left(\pi(q, x), v\right),
\]
or, in the other notation:
\[
q|_{xv}=q|_x|_v, \qquad q(xv)=q(x)q|_x(v).
\]

 We also put $q|_\emp=q$, $q(\emp)=\emp$.

Hence we get for every state $q$ a map $v\mapsto q(v)$, defining
the \emph{action of the state $q$} on the words. It is easy to see
that we have
\[
q_1q_2|_v=q_1|_{q_2(v)}g_2|_v,\qquad q(vw)=q(v)q|_v(w)
\]
for all $q, q_1, q_2\in Q$ and $v, w\in\xs$. Here $q_1q_2$ is the
product of transformations $q_1$ and $q_2$, i.e.,
$q_1q_2(w)=q_1\left(q_2(w)\right)$.

The above definitions imply the following description of
self-similar actions in terms of automata theory.

\begin{proposition}
\label{pr:full} A faithful action of a group $G$ on the set $\xs$
is self-similar if and only if there exists an automaton with the
set of states $G$ such that the action of the states of the
automaton on $\xs$ coincides with the original action of $G$.
\end{proposition}

The automaton from Proposition~\ref{pr:full} is called
\emph{complete automaton} of the action.

It is convenient to represent automata by their \emph{Moore
diagrams}. If $\mathcal{A}=\langle Q, \pi,\lambda\rangle$ is an
automaton, then its Moore diagram is a directed graph with the set
of vertices $Q$ in which we have for every pair $x\in X, q\in Q$
an arrow from $q\in Q$ into $\pi(q, x)$ labelled by the pair of
letters $(x; \lambda(q, x))$.

Let $q\in Q$ be a state and let $v\in\xs$ be a word. In order to
find the image $q(v)$ of the word $v$ under the action of the
state $q$ one needs to find a path in the Moore diagram, which
starts at the state $q$ with the consecutive labels of the form
$(x_1; y_1), (x_2; y_2), \ldots (x_n; y_n)$, where $x_1x_2\ldots
x_n=v$, then $q(v)=y_1y_2\ldots y_n$.

\begin{defi}
We say that an automaton $\mathcal{A}=\langle Q, \pi,
\lambda\rangle$ \emph{has finite nucleus} if there exists its
finite sub-automaton $\nuke\subset Q$ such that for every
$q\in\mathcal{A}$ there exists $n\in\mathbb{N}$ such that
$q|_v\in\nuke$ for all $v\in\xs$ such that $|v|\ge n$.

A self-similar action of a group $G$ on $\xs$ is said to be
\emph{contracting} if its full automaton has a finite nucleus.
\end{defi}

In general, if $\mathcal{A}$ is an automaton, then its
\emph{nucleus} is the set
\[
\nuke=\bigcup_{q\in
Q}\bigcap_{n\in\mathbb{N}}\left\{q|_v\;:\;v\in\xs, |v|\ge
n\right\}.
\]

For more on contracting actions, see the
papers~\cite{fractal_gr_sets,nekr:virt}.


\section{Limit spaces}

One of important properties of contracting actions is there strong
relation to Dynamical Systems, exhibited in the following notion
of \emph{limit space}.

Denote by $\xmo$ the set of all infinite to the left sequences of
the form $\ldots x_2x_1$, where $x_i$ are letters of the alphabet
$X$. We introduce on the set $\xmo$ the topology of the infinite
power of the discrete set $X$. Then the space $\xmo$ is a compact
totally disconnected metrizable topological space without isolated
points. Thus it is homeomorphic to the Cantor space.

\begin{defi}
Let $(G, \xs)$ be a contracting group action over the alphabet
$X$. We say that two points $\ldots x_2x_1, \ldots y_2y_1\in\xmo$
are \emph{asymptotically equivalent} (with respect to the action
of the group $G$) if there exists a bounded sequence
$\{g_k\}_{k\ge 1}$ of group elements such that for every
$k\in\mathbb{N}$ we have
\[
g_k(x_k\ldots x_1)=y_k\ldots y_1.
\]
\end{defi}

Here a sequence $\{g_k\}_{k\ge 1}$ is said to be \emph{bounded} if
the set of its values is finite.

It is easy to see that the defined relation is an equivalence. The
quotient of the space $\xmo$ by the asymptotic equivalence
relation is called the \emph{limit space} of the action and is
denoted $\lims$.

We have the following properties of the limit space
(see~\cite{nekr:limsp}).

\begin{theorem}
The asymptotic equivalence relation is closed and has finite
equivalence classes. The limit space $\lims$ is metrizable and
finite-dimensional. The shift $\sigma:\ldots x_2x_1\mapsto\ldots
x_3x_2$ induces a continuous surjective map $\si:\lims\arr\lims$.
\end{theorem}

We will also use the following description of the asymptotic
equivalence relation (for the proof see~\cite{nekr:limsp}).

\begin{proposition}
\label{pr:moorelim} Two sequences $\ldots x_2x_1, \ldots
ay_2y_1\in\xmo$ are asymptotically equivalent if and only if there
exists a sequence $g_1, g_2, \ldots$ of elements of the nucleus
such that $g_i\cdot x_i=y_i\cdot g_{i-1}$, i.e., if there exists a
left-infinite path $\ldots e_2, e_1$ in the Moore diagram of the
nucleus such that the edge $e_i$ is labelled by $(x_i; y_i)$.
\end{proposition}

A \emph{left-infinite path} in a directed graph is a sequence
$\ldots e_2, e_1$ of its arrows such that beginning of $e_i$ is
equal to the end of $e_{i+1}$. The end of the last edge $e_1$ is
called the \emph{end} of the left-infinite path.


The dynamical system $(\lims, \si)$ has a special Markov partition
coming from the its presentation as a shift-invariant quotient.

\begin{defi}
For every finite word $v\in\xs$ the respective \emph{tile}
$\til_v$ is the image of the cylindrical set $\xmo v$ in the limit
space $\lims$. We say that $\til_v$ is a tile of the level number
$|v|$.
\end{defi}

We have the following obvious properties of the tiles.

\begin{enumerate}
  \item Every tile $\til_v$ is a compact set.
  \item $\si\left(\til_{vx}\right)=\til_v$.
  \item $\til_{v}=\cup_{x\in X}\til_{xv}$.
\end{enumerate}

In particular, the image of a tile $\til_v$ of $n$th level under
the shift map $\si$ is a union of $d$ tiles $\til_u$ of the $n$th
level, i.e., that the tiles of one level for a Markov partition of
the dynamical system $(\lims, \til)$.

Actually, a usual definition of a Markov partition requires that
two tiles do not overlap, i.e., that they do not have common
interior points. We have the following criterion (for a proof see
also~\cite{nekr:limsp}).

We say that a self-similar action satisfies the \emph{open set
condition} if for every $g\in G$ there exists $v\in\xs$ such that
$g|_v=1$.

\begin{theorem}
\label{th:openset} If a contracting action of a group $G$ on $\xs$
satisfies the open set condition then for every $n\ge 0$ and for
every $v\in X^n$ the boundary of the tile $\til_v$ is equal to the
set
\[
\partial\til_v=\bigcup_{u\in X^n, u\ne v}\til_u\cap\til_v,
\]
and the tiles of one level have disjoint interiors.

If the action does not satisfy the open set condition, then there
exists $n\in\mathbb{N}$ and a tile of $n$th level, which is
covered by other tiles of $n$th level.
\end{theorem}

\section{Post-critically finite limit spaces}

Following~\cite{kigami:anal_fract}, we adopt the following
definition.

\begin{defi}
We say that a contracting action $(G, \xs)$, has a
\emph{post-critically finite (p.c.f.) limit space} if intersection
of every two different tiles of one level is finite.
\end{defi}

We obtain directly from Theorem~\ref{th:openset} that a
contracting action has a p.c.f.\ limit space if and only if it
satisfies the open set condition and the boundary of every tile is
finite.

The following is an easy corollary of Theorem~\ref{th:openset} and
Proposition~\ref{pr:moorelim}.

\begin{proposition}
\label{pr:intersection} The image of a sequence $\ldots
x_{n+1}x_n\ldots x_1\in\xmo$ belongs to the boundary of the tile
$\til_{x_n\ldots x_1}$ if and only if there exists a sequence
$\{g_k\}$ of elements of the nucleus such that $g_{k+1}\cdot
x_{k+1}=x_k\cdot g_k$ and $g_n(x_n\ldots x_1)\ne x_n\ldots x_1$.
\end{proposition}

This gives us an alternative way of defining p.c.f.\ limit spaces.

\begin{corollary}
\label{cor:pcflim} A contracting action $(G, \xs)$ has a p.c.f.\
limit space if and only if there exists only a finite number of
left-infinite paths in the Moore diagram of its nucleus which end
in a non-trivial state.
\end{corollary}

\begin{proof}
We say that a sequence $\ldots x_2x_1\in\xmo$ is read on a
left-infinite path $\ldots e_2e_1$, if the label of the edge $e_i$
is $(x_i;y_i)$ for some $y_i\in X$. If the path $\ldots e_2e_1$
passes through the states $\ldots g_2g_1g_0$ (here $g_i$ is the
beginning and $g_{i-1}$ is the end of the edge $e_i$), then the
state $g_{n-1}$ is uniquely defined by $g_n$ and $x_n$, since
$g_{n-1}=g_n|_{x_n}$. Consequently, any given sequence $\ldots
x_2x_1$ is read not more than on $|\nuke|$ left-infinite paths of
the nucleus $\nuke$. In particular, every asymptotic equivalence
class on $\xmo$ has not more than $|\nuke|$ elements.

For every non-trivial state $g\in\nuke$ denote by $B_g$ be the set
of sequences, which are read on the left-infinite paths of the
nucleus, which end in $g$.

Suppose that there is infinitely many left-infinite paths in the
nucleus ending in a nontrivial state. Then there exists a state
$g\in\nuke\setminus\{1\}$ for which the set $B_g$ is infinite.

Since the state $g$ is non-trivial, there exists a word $v\in\xs$
such that $g(v)\ne v$. Then for every $\ldots x_2x_1\in B_g$,
there exists a sequence $\ldots y_2y_1$ such that $\ldots x_2x_1v$
is asymptotically equivalent to $\ldots y_2y_1g(v)$. Hence, every
point of $\lims$ represented by a sequence from $B_gv$ belongs
both to $\til_v$ and to $\til_{g(v)}$. This show that the
intersection $\til_v\cap\til_{g(v)}$ is infinite, since the
asymptotic equivalence classes are finite.

On the other hand, Proposition~\ref{pr:intersection} shows, that
if the sequence $\ldots x_2x_1v$, represents a point of the
intersection $\til_v\cap\til_u$ for $u\in X^{|v|}, u\ne v$, then
the sequence $\ldots x_2x_1$ is read on some path of the nucleus,
which ends in a non-trivial state. Therefore, if the intersection
$\til_v\cap\til_u$ is infinite, then the set of left-infinite
paths in the nucleus is infinite.
\end{proof}

\begin{proposition}
\label{pr:1dim} If the limit space of a contracting action is
post-critically finite, then its topological dimension is $\le 1$.
\end{proposition}

\begin{proof}
We have to prove that every point $\zeta\in\lims$ has a basis of
neighborhoods with 0-dimensional boundaries.

Let $T_n(\zeta)$ be the union of the tiles of $n$th level,
containing $\zeta$. It is easy to see that
$\{T_n(\zeta)\;:\;n\in\mathbb{N}\}$ is a base of neighborhoods of
$\zeta$.
\end{proof}

\section{Automata with bounded cyclic structure}

We take Corollary~\ref{cor:pcflim} as a justification of the
following definition.

\begin{defi}
A self-similar contracting group is said to be
\emph{post-critically finite} (\pcf\ for short) if there exists
only a finite number of inverse paths in the nucleus ending at a
non-trivial state.
\end{defi}

A more precise description of the structure of the nucleus of a
\pcf\ group is given in the next proposition.

Recall, that an automatic transformation $q$ of $\xs$ is said to
be \emph{finitary} (see~\cite{grineksu_en}) if there exists
$n\in\mathbb{N}$ such that $q|_v=1$ for all $v\in X^n$ (then
$q(x_1\ldots x_m)=q(x_1\ldots x_n)x_{n+1}\ldots x_m$). The minimal
number $n$ is called \emph{depth} of $q$.

The set of all finitary automatic transformations of $\xs$ is a
locally finite group. If $G$ is a finite subgroup of the group of
finitary transformations, then the \emph{depth} of $G$ is the
greatest depth of its elements.

If we have a subset $A$ of the vertex set of a graph $\Gamma$,
then we consider it to be a subgraph of $\Gamma$, taking all the
edges, which start and end at the vertices of $A$.

We say that a directed graph is a \emph{simple cycle} if its
vertices $g_1, g_2, \ldots, g_n$ and edges $e_1, e_2, \ldots, e_n$
can be indexed so that $e_i$ starts at $g_i$ and ends at $g_{i+1}$
(here all $g_i$ and all $e_i$ are pairwise different and
$g_{n+1}=g_1$).

\begin{proposition}
\label{pr:nuclstruct} Let $\nuke$ be the nucleus of a \pcf\ group,
and let $\nuke_0$ be the subgraph of finitary elements of $\nuke$
and $\nuke_1=\nuke\setminus\nuke_0$. Then $\nuke_1$ is a disjoint
union of simple cycles.
\end{proposition}

\begin{proof}
The set $\nuke_0$ is obviously a sub-automaton, i.e., for every
$g\in\nuke_0$ and $x\in X$ we have $g|_x\in\nuke_0$. It follows
then from the definition of a nucleus that every vertex of the
graph $\nuke_1$ has an incoming arrow. This means that every
vertex of the graph $\nuke_1$ is an end of a left-infinite path.
On the other hand, there exists for every $g\in\nuke_1$ at least
one $x\in X$ such that $g|_x\in\nuke_1$, since all elements of
$\nuke_1$ are not finitary. Thus, every vertex of $\nuke_1$ has an
outgoing arrow and is a beginning of a right-infinite path.

Let $g$ be an arbitrary vertex of the graph $\nuke_1$. We have a
left-infinite path $\gamma_-$ ending in $g$. Suppose that we have
a (pre-)periodic right-infinite path $\gamma_+$ starting at $g$,
i.e., a path of the form $\gamma_+=qppp\ldots=qp^\omega$, were $q$
is a finite path, $p$ is a finite simple cycle and the set of
edges of the paths $p$ and $q$ are disjoint. Note that there
always exists a (pre-)periodic path beginning at $g$.

If $q$ is not empty, then we get an infinite set of different
left-infinite paths in the graph $\nuke_1$:
$\{\gamma_-qp^n\}_{n\in\mathbb{N}}$, what contradicts to the
post-critical finiteness of the action.

Hence the pre-period $q$ is empty. In particular, every element of
$\nuke_1$ belongs to a finite cycle, i.e., for every $g\in\nuke_1$
there exists $v\in\xs$ such that $g|_v=g$.

Suppose now that there exist two different letters $x, y\in X$
such that $g|_x$ and $g|_y$ belong to $\nuke_1$. The element
$g|_x$ belongs to a finite cycle $p_x$ in $\nuke_1$. The cycle
$p_x$ must contain the element $g$, otherwise we get a strictly
pre-periodic path starting at $g$. Similarly, there exists a cycle
$p_y$ , which contains $g$ and $g|_y$. The cycles $p_x$ and $p_y$
are different and intersect in the vertex $g$. Hence we get an
infinite set of left-infinite paths in $\nuke_1$ of the form
$\ldots p_3p_2p_1$, where $p_i$ are either $p_x$ or $p_y$ (seen as
paths starting at $g$) in an arbitrary way.

Hence, for every $g\in\nuke_1$ there exists only one letter $x\in
X$ such that $g|_x\in\nuke_1$. This (together with the condition
that every vertex of $\nuke_1$ has an incoming edge) implies that
$\nuke_1$ is a disjoint union of simple cycles.
\end{proof}

The following notion was defined and studied by Said Sidki
in~\cite{sid:cycl}.

\begin{defi}
We say that an automatic transformation $q$ is \emph{bounded} if
the sequence $\theta(k, q)$ is bounded, where $\theta(k, q)$ is
the number of words $v\in X^k$ such that $q|_v$ acts non-trivially
on the first level $X^1$ of the tree $\xs$.
\end{defi}

The following proposition is proved in~\cite{sid:cycl}
(Corollary~14).

\begin{proposition}
An automatic transformation is bounded if and only if it is
defined by a finite automaton in which every two non-trivial
cycles are disjoint and not connected by a directed path.
\end{proposition}

Here a cycle is \emph{trivial} if its only vertex is the trivial
state. In particular, every finitary transformation is bounded,
since it has no non-trivial cycles.

\begin{theorem}
\label{th:main} The set $\B$ of all bounded automorphisms of the
tree $\xs$ is a group.

A finitely generated self-similar automorphism group $G$ of the
tree $\xs$ has a p.c.f.\ limit space if and only if it is a
subgroup of $\B$. In particular, every finitely generated
self-similar subgroup of $\B$ is contracting.
\end{theorem}

\begin{proof}
The fact that $\B$ is a group, is proved in~\cite{sid:cycl}. We
have also proved that the nucleus of every p.c.f.\ group $G$ is a
subset of $\B$. This implies that $G$ is a subgroup of $\B$.

In the other direction, suppose that we have a self-similar
finitely generated subgroup $G\le\B$. Then $G$ is generated by a
finite automaton $S$ whose all non-trivial cycles are disjoint.
Let $S_0$ be the subautomaton of all finitary tranformations, and
let $S_2=S\setminus S_0$. Then all non-trivial cycles belong to
$S_2$. Let $S_1$ be the union of all these cycles.

Let $g\in S_1$ and $v\in\xs$ be arbitrary. Then either $g|_v$
belongs to the same cycle as $g$, or $g|_v\notin S_1$, since no
two different cycles of $S_1$ can be connected by a directed path.
If $g|_v\notin S_1$, then all states $g|_{vu}$ of $g|_v$ do not
belong to $S_1$. But this is possible only when $g|_v\in S_0$.
Therefore, there exists $m\in\mathbb{N}$ such that for every $g\in
S$ and every $v\in X^m$ either $g|_v\in S_0$, or $g|_v\in S_1$.
Then the group $G_1=\langle G|_{X^m}\rangle$ is also self-similar
and is generated by a subset of the set $S_0\cup S_1$. The group
$G$ is contracting if and only if $G_1$ is contracting. Their
nuclei will coincide. Therefore, if we prove our theorem for
$G_1$, then it will follow for $G$, so we assume that $S_2=S_1$.

Let $n_1$ be the least common multiple of the lengths of cycles in
$S_1$. Then for every $u\in X^{n_1}$ and $s\in S_1$ we have either
$s|_u\in S_0$ or $s|_u=s$. Moreover, it follows from the
conditions of the theorem that the word $u\in X^{n_1}$ such that
$s|_u=s$ is unique for every $s\in S_1$.

Let $\nuke_1$ be the set of all elements $h\in G\setminus 1$ such
that there exists one word $u(h)\in X^{n_1}$ such that
$h|_{u(h)}=h$ and for all the other words $u\in X^{n_1}$ the
restriction $h|_u$ belongs to $\langle S_0\rangle$. It is easy to
see that the set $\nuke_1$ is finite (every its element $h$ is
uniquely defined by the permutation it induces on $X^{n_1}$ and
its restrictions in the words $u\in X^{n_1}$, note also that the
group $\langle S_0\rangle$ is finite).

Let us denote by $l_1(g)$ the minimal number of elements of
$S_1\cup S_1^{-1}$ in a decomposition of $g$ into a product of
elements of $S\cup S^{-1}$.

Let us prove that there exists for every $g\in G$ a number $k$
such that for every $v\in X^{n_1k}$ the restriction $g|_v$ belongs
to $\nuke_1\cup\{S_0\}$. We will prove this by induction on
$l_1(g)$.

If $l_1(g)=1$, then $g=h_1sh_2$, where $h_1, h_2\in\langle
S_0\rangle$ and $s\in S_1$. The elements $h_1, h_2$ are finitary,
thus there exists $k$ such that for every $v\in X^{n_1k}$ the
restriction $h_i|_v$ is trivial. Then we have
$h_1sh_2|_v=s|_{h_2(v)}$, thus $g|_v$ is either equal to
$s\in\nuke_1$ or belongs to $S_0\cup S_0^{-1}$. Thus the claim is
proved for the case $l_1(g)=1$.

Suppose that the claim is proved for all elements $g\in G$ such
that $l_1(g)<m$. Let $g=s_1s_2\ldots s_k$, where $s_i\in S\cup
S^{-1}$. For every $u\in X^{n_1}$ the restriction $s_i|_u$ is
equal either to $s_i$ or belongs to $S_0$. Consequently, either
$g|_u=g$ for one $u$ and $g|_u\in\langle S_0\rangle$ for all the
other $u\in X^{n_1}$, or $l_1(g|_u)<l_1(g)$ for every $u\in
X^{n_1}$. In the first case we have $g\in\nuke_1$ and in the
second we apply the induction hypothesis, and the claim is proved.

Consequently, the group $G$ is contracting with the nucleus equal
to a subset of the set $\{g|_v\;:\;g\in\nuke_1, v\in\xs,
|v|<n_1\}$. Note that any restriction $g|_v$ of an element of
$\nuke_1$ either belongs to $\nuke_1$ or is finitary.

Let us prove that the limit space of the group $G$ is \pcf.
Suppose that we have a left-infinite path in the nucleus of the
group. Let
\[
\ldots h_3, h_2, h_1
\]
be the elements of the nucleus $\nuke$ it passes through and let
the letters $\ldots, x_3, x_2, x_1$ be the letters labeling its
edges. In other words, we have
\[
h_n=h_{n+1}|_{x_n}
\]
for every $n\ge 1$.

The number of possibilities for $h_n$ is finite, thus it follows
from the arguments above that every element $h_i$ belongs to
$\nuke_1\cup\langle S_0\rangle$. The elements of $\langle
S_0\rangle$ can belong only to the ending of the sequence $h_i$ of
the length not greater than the depth of the group $\langle
S_0\rangle$. The rest of the sequences $h_i$ and $x_i$ is periodic
with period $n_1$. Hence, there exists only a finite number of
possibilities for such a sequence, and the limit space of the
group is \pcf
\end{proof}

\begin{corollary}
The word problem is solvable in polynomial time for every finitely
generated subgroup of $\mathcal{B}$.
\end{corollary}

\begin{proof}
The word problem in every finitely generated contracting group is
solvable in polynomial time (see~\cite{nekr:virt}).
\end{proof}

\bibliographystyle{plain}
\def\cprime{$'$} \def\cprime{$'$}

\end{document}